\documentclass[10pt]{article}
\usepackage{amssymb,amsmath,amsthm}
\newcommand{\R}{\mathbb{R}}
\begin{document}
\title{\Large\bf{Existence of solutions for Kirchhoff-type fractional Dirichlet problem  with $p$-Laplacian}
 }
\date{}
\author {  Danyang Kang$^{1}$,  Cuiling Liu$^{1}$, \ Xingyong Zhang$^{1,2}$\footnote{Corresponding author, E-mail address: zhangxingyong1@163.com}\\\
        {\footnotesize $^1$Faculty of Science, Kunming University of Science and Technology, Kunming, Yunnan, 650500, P.R. China.}\\
        {\footnotesize  $^2$School of Mathematics and Statistics, Central South University,
        Changsha, Hunan, 410083, P.R. China.}\\}

 \date{}
 \maketitle

 \begin{center}
 \begin{minipage}{15cm}
 \par
 \small  {\bf Abstract:} In this paper, we investigate the existence  of solutions for a class of $p$-Laplacian fractional order Kirchhoff-type system with  Riemann-Liouville fractional derivatives and a parameter $\lambda$.
 By mountain pass theorem, we obtain that  system has at least one non-trivial weak solution $u_\lambda$ under some local superquadratic conditions for each given  large parameter $\lambda$. We get a
concrete lower bound of the parameter $\lambda$, and then obtain two estimates of weak solutions $u_\lambda$. We also obtain that $u_\lambda\to 0$ if $\lambda$ tends to $\infty$. Finally, we present an example as an application of  our  results.
 \end{minipage}
 \end{center}
  \allowdisplaybreaks
 \vskip2mm
 {\section{Introduction and main results }}
\par
In this paper, we are concerned with the following system
\begin{equation}\label{aa2}
 \begin{cases}
 A(u(t))[_{t}D_{T}^{\alpha}\phi_{p}(_{0}D_{t}^{\alpha}u(t))+V(t)\phi_{p}(u(t))]
     =\lambda\nabla F(t,u(t)), \;\;\; \mbox{a.e. } t\in[0,T],\\
   u(0)=u(T)=0,
   \end{cases}
 \end{equation}
 where
 $$
 A(u(t))= \left[a+b\int_{0}^{T}(|_{0}D_{t}^{\alpha}u(t)|^{p}
     +V(t)|u(t)|^{p})dt\right]^{p-1},
 $$
  $a,b,\lambda>0$, $p>1$ and $1/p<\alpha\leq1$ are constants, $p$ is an integer, $u(t)=(u_1(t),\cdots, u_N(t))^\tau\in \R^N$ for a.e. $t\in [0,T]$, $T>0$, and $N$ is a given positive integer, $(\cdot)^\tau$ denote the transpose of a vector, $V(t)\in C([0,T],\R)$ with $\min\limits_{t\in[0,T]}V(t)>0$, $_{0}D_{t}^{\alpha}$ and $_{t}D_{T}^{\alpha}$ are the left and right Riemann-Liouville fractional derivatives, respectively, $\phi_{p}(s):=|s|^{p-2}s$, $\nabla F(t,x)$ is the gradient of $F$ with respect to $x=(x_1,\cdots, x_N)\in \R^N$, that is, $\nabla F(t,x)=(\frac{\partial F}{\partial x_1},\cdots,\frac{\partial F}{\partial x_N})^\tau$,  and $F:[0,T]\times\R^N\rightarrow \R$ satisfies the following condition:
\vskip2mm
\noindent
 {\it (H0) there exists a constant $\delta>0$ such that $F(t,x)$ is continuously differentiable in $x\in \R^N$ with $|x|\le \delta$ for a.e. $t \in [0,T]$, measurable
in $t$ for every $x\in \R^N$ with $|x|\le \delta$, and there exist $a\in C(\R^+,\R^+)$ and
$b\in L^1([0,T];\R^+)$ such that
$$
|F(t, x)|, |\nabla F(t,x)|\le a(|x|)b(t)
$$
for all $x\in \R^N$ with $|x|\le \delta$ and a.e. $t \in [0,T]$.
}
\par
 When $\alpha=1$, the  operator $_{t}D_{T}^{\alpha}(_{0}D_{t}^{\alpha}u(t))$ reduces to
the usual second order differential operator $-d^2/dt^2$. Hence, if $\alpha=1$,  $p=2$, $N=1$, $\lambda=1$ and $V(t)=0$ for a.e. $t\in [0,T]$, system (\ref{aa2}) becomes the equation with Dirichlet  boundary condition
\begin{equation}\label{aa6}
 \begin{cases}
  -\left(a+b\int_{0}^{T}|u'(t)|^{2}
     dt\right)u''(t)
     =f(t,u(t)), \;\;\; \mbox{a.e. } t\in[0,T],\\
   u(0)=u(T)=0,
   \end{cases}
 \end{equation}
 where $f(t,x)=\frac{\partial F(t,x)}{\partial x}$ and $F:[0,T]\times\R\to \R$. It is well known that equation (\ref{aa6}) is related to  the stationary problem of a classical model introduced by Kirchhoff \cite{Kir}. To be precise, in  \cite{Kir}, Kirchhoff introduced the  model
 \begin{equation}\label{aa7}
 \rho \frac{\partial^2 u}{\partial t^2}=\left(P_0+\frac{Eh}{2L}\int_0^L\left(\frac{\partial u}{\partial y}\right)^2dy\right)\frac{\partial^2 u}{\partial y^2},
 \end{equation}
where $0\le y \le L$, $t \ge 0$, $u$ is the lateral deflection, $\rho$ is the mass density, $h$ is the cross-sectional area, $L$ is the length, $E$ is the Young¡¯s modulus and $P_0$ is the initial axial tension. ({\bf Notations}: in model (\ref{aa7}), (\ref{kk4}) and (\ref{kk5}) below, $t$ is time variable and $y$ is
spatial variable, which are conventional notations in partial differential equations. One need to distinguish them to $t$ in (\ref{aa2}), (\ref{aa6}),  (\ref{kk1}), (\ref{kk2}) and (\ref{kk3}) below, where $t$ corresponds to the spatial variable $x$). The model (\ref{aa7}) is used to describe small vibrations of an elastic stretched string.
Equation (\ref{aa7}) has been studied extensively, for instance, \cite{Ber}, \cite{ANC}, \cite{GRE}, \cite{POH}, \cite{Arosio A1996}, \cite{Cava}, \cite{ANC2}, \cite{Alves2005}, \cite{Perera2006},  \cite{He2014}, \cite{Anran2015}, \cite{Tang2016}, \cite{Yang2019} and reference therein. For $p>1$, the reader can consult \cite{Dreher2006}, \cite{Francisco2006}, \cite{Liu2010}, \cite{Pucci2010}, \cite{Liu2013} and references therein.

\par
When $\alpha<1$,  $_{0}D_{t}^{\alpha}$ and $_{t}D_{T}^{\alpha}$ are the left and right Riemann-Liouville fractional derivatives, respectively, which has been given some physical interpretations in \cite{Heymans2006}. Moreover, it is also applied to describe the
 anomalous diffusion, L\'evy flights and traps in \cite{Ervin2007} and \cite{Zhuang2008}. Fractional differential equations have been proved to provide a natural framework  in the modeling of many real phenomena such as viscoelasticity, neurons, electrochemistry, control, porous media, electromagnetic (the reader can consult \cite{Jiao F2012} in which a collection of references is given).  In \cite{Jiao F2012}, Jiao and Zhou considered the  system
 \begin{equation}\label{kk1}
 \begin{cases}
_{t}D_{T}^{\alpha}(_{0}D_{t}^{\alpha}u(t)) =\nabla F(t,u(t)), \;\;\; \mbox{a.e. } t\in[0,T],\\
   u(0)=u(T)=0.
   \end{cases}
 \end{equation}
They successfully applied critical point theory to investigate the existence of weak solutions for system (\ref{kk1}). To be precise, they obtained that system (\ref{kk1}) has at least one weak solution when $F$ has a quadratic growth or a superquadratic growth by using the least action principle and mountain pass theorem. Subsequently, this topic related to system (\ref{kk1}) attracted lots of attention, for example, \cite{Zhou2014}, \cite{Zhao2015}, \cite{Bonanno2014}, \cite{Zhao2017}, \cite{Xie2018}, \cite{Nyamoradia2017} and references therein. It is obvious that system (\ref{aa2}) is much more complicated than  system (\ref{kk1}) since the appearance of nonlocal term $A(u(t))$ and $p$-Laplacian term $\phi_p(s)$. Recently, in \cite{Chai2018}, the following fractional Kirchhoff equation with Dirichlet boundary condition was investigated
\begin{equation}\label{kk2}
 \begin{cases}
\left(a+b\int_{0}^{T}|_{0}D_{t}^{\alpha}u(t)|^{2}dt\right){_{t}D}_{T}^{\alpha}(_{0}D_{t}^{\alpha}u(t))+\lambda V(t)u(t)
     =f(t,u(t)), \;\;\; \mbox{a.e. }t\in[0,T],\\
   u(0)=u(T)=0,
   \end{cases}
 \end{equation}
where $a,b,\lambda>0$, $f\in C([0,T]\times\R,\R)$. By using the mountain pass theorem in \cite{Eke} and
the linking theorem in \cite{Li}, the authors established some existence results of nontrivial solutions for system (\ref{kk2}) if
$f$ satisfies\\
{\it (f1) there exist constants $\mu> 4$, $0 < \tau < 2$ and a nonnegative function $g \in L^{\frac{2}{2-\tau}}$ such that
$$
F(t, x)- \frac{1}{\mu}f(t, x)x\le g(t)|x|^\tau, \ \ \ \mbox{for a.e.  } t \in [0,T], x \in\R;
$$
(f2) there exists $\theta > 2$ such that $\lim_{|x|\to\infty} \inf_{t\in[0,T]}
\frac{F(t,x)}{|x|^\theta} > 0$;\\
(or
(f2)$'$ there exists $\theta> 4$ such that $\lim_{|x|\to\infty} \inf_{t\in[0,T]}
\frac{F(t,x)}{|x|^\theta} > 0$);\\
(f3) there exists $\sigma  > 2$ such that  $\lim_{|x|\to 0} \sup_{t\in[0,T]}
\frac{F(t,x)}{|x|^\sigma}<\infty,$\\
}
and some other reasonable conditions.
\par
In \cite{Chen2018}, Chen-Liu investigated the Kirchhoff-type
fractional Dirichlet problem with $p$-Laplacian
\begin{equation}\label{kk3}
\begin{cases}
\left(a+b\int_{0}^{T}|_{0}D_{t}^{\alpha}u(t)|^{p}dt\right)^{p-1}{_{t}D}_{T}^{\alpha}\phi_p(_{0}D_{t}^{\alpha}u(t))
     =f(t,u(t)), \;\;\; t\in(0,T),\\
   u(0)=u(T)=0.
   \end{cases}
 \end{equation}
where $a,b,\lambda>0$, $f\in C^1([0,T]\times\R,\R)$. By using the Nehari method, they established the existence result of ground state solution for system (\ref{kk3}) if $f$ satisfies\\
{\it  (f4) $f (t, x) = o(|x|^{p-1})$ as $|x|\to 0$ uniformly for all $ t \in [0,T]$,}\\
and the following well-known Ambrosetti-Rabinowitz (AR for short) condition\\
{\it (AR) there exist two constants $\mu>p^2$, $R>0$ such that
$$
0 < \mu F(t, x) \le  xf (t, x), \ \ \mbox{for }\forall t \in [0,T], x\in \R \mbox{ with } |x|\ge R,
$$
where $F(t, x) =\int_0^x f(t,s)ds$}, and some additional conditions. It is easy to see that all of these conditions (f1), (f2), (f2)$'$ and (AR) imply that  $F(t, x)$ needs to have a growth near the infinity for $x$, and (f3) and (f4) imply that  $F(t, x)$ needs to have a growth near $0$ for $x$.
\par
 In this paper, we investigate the existence and concentration of solutions for system (\ref{aa2}) under local assumptions only near $0$  for the nonlinear term $F$. Our work is mainly motivated by \cite{costa} and \cite{Anran2015}. In \cite{costa}, Costa and Wang investigated the multiplicity of both signed and sign-changing solutions for the one-parameter family of elliptic problems
 \begin{equation}\label{kk4}
\begin{cases}
\Delta u = \lambda f(u) \ \ \mbox{ in } \Omega\\
  u(y) = 0 \ \ \mbox{ in }\partial \Omega,
   \end{cases}
 \end{equation}
 where $\lambda>0$ is a parameter, $\Omega$ is a bounded smooth domain in $\R^N (N \ge 3)$ and $f\in C^1(\R,\R)$. They assumed that the nonlinearity $f(u)$ has superlinear growth in a neighborhood of $u = 0$ and then obtained the number of signed and sign-changing solutions which are dependent on the parameter $\lambda$. The idea in \cite{costa} has been applied to some different problems, for example, \cite{Medeiros2007} and \cite{Nikolaos2018} for  quasilinear elliptic problems with $p$-Laplacian operator, \cite{XU2016} for an  elliptic problem with fractional Laplacian operator, \cite{Joao2008} for Schr$\ddot{\mbox{o}}$dinger equations, \cite{He2014} for Neumann problem with nonhomogeneous differential operator and critical growth, and \cite{Huang2019} for quasilinear Schr$\ddot{\mbox{o}}$dinger equations. Especially, in \cite{Anran2015}, Li and Su investigated the Kirchhoff-type equations
 \begin{equation}\label{kk5}
\begin{cases}
-\left[1+\int_{\R^3}(|\nabla u|^{2}+V(y)u^2)dy\right][\Delta u+V(y)u]
     =\lambda Q(y)f(u), \ \ y\in\R^3,\\
   u(y)\to 0, \ \ \mbox{as }|y|\to \infty,
   \end{cases}
 \end{equation}
where $\lambda>0$, $V,Q$ are radial functions and $f\in C((-\delta_0,\delta_0),\R)$ for some $\delta_0>0$. Via the idea in \cite{costa}, they also established the existence result of solutions when $f(u)$ has superlinear growth in a neighborhood of $u = 0$.
 It is worthy to note that $\lambda$ usually needs to be sufficiently large, that is, $\lambda$ has a lower bound $\lambda^*$. However, the concrete values of $\lambda^*$ are not given in these references.
 Similar to system (\ref{kk5}), comparing with equation (\ref{kk2}) and equation (\ref{kk3}), we add a nonlocal term $\int_0^T V(t)|u(t)|^p dt$ in  system (\ref{aa2}) where $\min_{t\in [0,T]}V(t)>0$, and multiply $V(t)\phi_p(u(t))$ by the nonlocal part $A(u(t))$. Moreover, we  consider the high-dimensional case, that is, $N\ge 1$. Since $\min_{t\in [0,T]}V(t)>0$, system (\ref{aa2}) is different from  equation (\ref{aa6}),  (\ref{kk2}), (\ref{kk3}) and system  (\ref{kk1}). More importantly, we present a concrete value of the lower bound $\lambda^*$ for system (\ref{aa2}) and then obtain two estimates of the solutions family $\{u_\lambda\}$ for all $\lambda>\lambda^*$.
Next, we make some assumptions for $F$.\\
{\it (H1)  there exist constants $q_{1}>p^2$, $q_{2}\in(p^2,q_{1})$, $M_1>0$ and $M_2>0$ such that
$$
 M_1 |x|^{q_1}\le F(t,x)\le M_2 |x|^{q_2}
$$
for all  $x\in \R^N$  with $|x|\le \delta$ and a.e. $t\in [0,T]$;\\
(H2) there exists a constant $\beta>p^{2}$  such that
$$
0\leq\beta F(t,x)\leq(\nabla F(t,x),x)
$$
for all  $x\in \R^N$  with $|x|\le \delta$ and a.e. $t\in [0,T]$.}
 \vskip0mm
 \noindent
 {\bf Theorem 1.1.} {\it Suppose that (H0)-(H2) hold. Then  system (\ref{aa2})  has at least a nontrivial weak solution $u_{\lambda}$ for all $\lambda>\lambda^*:=\max\{\Lambda_1,\Lambda_2,\Lambda_3\}$ and
\begin{eqnarray*}
&  &         \|u_{\lambda}\|_V^{p}
           \leq  \dfrac{p^{2}\theta}{a^{p-1}(\theta-p^{2})}  C_* \lambda^{-\frac{p-1}{q_{1}-p}}\leq \dfrac{p^{2}\theta}{a^{p-1}(\theta-p^{2})}  C_* \max\{\Lambda_1,\Lambda_2,\Lambda_3\}^{-\frac{p-1}{q_{1}-p}},\\
&  &            \|u_{\lambda}\|_\infty
           \leq  \frac{T^{\alpha-\frac{1}{p}}}{\Gamma(\alpha)(\alpha q-q+1)^{\frac{1}{q}}} \cdot \dfrac{p^{2}\theta}{a^{p-1}(\theta-p^{2})}  C_* \lambda^{-\frac{p-1}{q_{1}-p}}\\
&  &     \quad\quad\quad\quad     \leq \frac{T^{\alpha-\frac{1}{p}}}{\Gamma(\alpha)(\alpha q-q+1)^{\frac{1}{q}}} \cdot \dfrac{p^{2}\theta}{a^{p-1}(\theta-p^{2})}  C_* \max\{\Lambda_1,\Lambda_2,\Lambda_3\}^{-\frac{p-1}{q_{1}-p}},\\
&  & \lim_{\lambda\to \infty} \|u_{\lambda}\|_V=0=\lim_{\lambda\to \infty} \|u_{\lambda}\|_\infty,
 \end{eqnarray*}
 where $\theta=\min\{\beta,q_2\}$,} $q=\frac{p}{p-1}$,
 \begin{eqnarray}
 & & \label{P1}\|u_\lambda\|_{V}=\left(\int_0^T |_0D_t^\alpha u_\lambda(t)|^pdt+\int_0^T V(t)|u_\lambda(t)|^pdt\right)^{1/p},\ \ \|u_\lambda\|_{\infty}=\max_{t\in[0,T]}u_\lambda(t),\\
 & & \label{A1}\Lambda_1=\max\left\{\frac{V_\infty a^{p-1}(\Gamma(\alpha)(\alpha q-q+1)^{\frac{1}{q}}G_0)^{q_2-p}}{2p^2M_2T^{(\alpha-\frac{1}{p})(q_2-p)}(\delta \min\{1,V_\infty\} D)^{q_2-p}}, \frac{\dfrac{1}{bp^{2}}(a+\frac{b\delta^p}{G_0^p}\max\{1,V^\infty\}(D^{p}+G^{p}))^{p}}{M_1 \frac{\delta^{q_1}}{G_0^{q_1}}T^{1-\frac{q_1}{p}}D^{q_1} }\right\},\\
 & & \label{A2}\Lambda_2 =\left[a+b[\max\{1,V^\infty\}]^p\frac{\delta^{p}}{G_0^{p}}(D^p+G^p)\right]^{q_1(p-1)},\\
 & & \label{A3}\Lambda_3= \left(\frac{T^{p\alpha-1}}{\left[\Gamma(\alpha)(\alpha q-q+1)^{\frac{1}{q}}\right]^p} \cdot \dfrac{p^{2}\theta C_* }{a^{p-1}(\theta-p^{2})} \cdot \frac{2^p}{\delta^p} \right)^{\frac{q_1-p}{p-1}},\\
 & & V^\infty=\max_{t\in [0,T]} V(t),\ \ V_\infty=\min_{t\in [0,T]} V(t),\nonumber\\
 & & \label{kkk3}C_*= \left(\frac{1}{p(M_1q_{1})^{\frac{p}{q_{1}-p}}}-\frac{M_1}{(M_1q_{1})^{\frac{q_{1}}{q_{1}-p}}}\right)
       \left(\dfrac{[\max\{1,V^\infty\}]^{1/p}(D^p+G^p)^{1/p}}{T^{\frac{1}{q_1}-\frac{1}{p}}D}\right)^{\frac{pq_{1}}{q_{1}-p}} ,\\
 & & D= \begin{cases}
               \left(\frac{T^{p+1}}{\pi^{p+1}}\cdot\frac{2(p-1)!!}{p!!} \right)^{\frac{1}{p}},&  \mbox{\it if } p \mbox{ \it is odd},\\
               \left(\frac{T^{p+1}}{\pi^{p}}\cdot\frac{(p-1)!!}{p!!}\right)^{\frac{1}{p}}, & \mbox{\it if } p\  \mbox{\it is even},
               \end{cases}\nonumber\\
 &  & G=\bigg(\frac{T^{p+1-p\alpha}}{[\Gamma(2-\alpha)]^p(p+1-p\alpha)}\bigg)^{1/p},\nonumber\\
 &  & G_0=\frac{T^{\alpha-\frac{1}{p}}}{\Gamma(\alpha)(\alpha q-q+1)^{\frac{1}{q}}} G ,\nonumber\\
 &  & \Gamma(z)=\int_0^\infty t^{z-1}e^{-t}dt\ \ (\mbox{ for all } z>0).\nonumber
 \end{eqnarray}
\par
We organize this paper as follows. In section 2, we recall some preliminary results including the definitions of  Riemann-Liouville fractional derivatives and working spaces,  some conclusions for the working spaces and mountain pass theorem. In section 3, we give the proof of Theorem 1.1. In section 4, we apply Theorem 1.1 to an example and compute the value of lower bound $\lambda^*$.

 \vskip2mm
 {\section{Preliminaries}}
 In this section, we mainly recall some basic definitions and results.
 \vskip2mm
 \noindent
{\bf Definition 2.1.} (Left and Right Riemann-Liouville Fractional Integrals \cite{Kilbas A A A2006,Zhou2014}) Let $f$ be a function defined on $[a,b]$. The left and right Riemann-Liouville fractional integrals of order $\gamma>0$ for function $f$ denoted by $_aD_t^{-\gamma} f(t)$ and $_tD_b^{-\gamma} f(t)$ , respectively, are defined by
 \begin{eqnarray*}
  &  & _aD_t^{-\gamma} f(t)=\frac{1}{\Gamma(\gamma)}\int_a^t (t-s)^{\gamma-1} f(s)ds , \ t\in[a,b],\gamma>0,\\
  &  & _tD_b^{-\gamma} f(t)=\frac{1}{\Gamma(\gamma)}\int_t^b (s-t)^{\gamma-1} f(s)ds , \ t\in[a,b],\gamma>0,
 \end{eqnarray*}
  provided the right-hand sides are pointwise defined on $[a,b]$, where $\Gamma>0$ is the Gamma function.
  \vskip2mm
 \noindent
 {\bf Definition 2.2.} (Left and Right Riemann-Liouville Fractional Derivatives\cite{Kilbas A A A2006, Zhou2014}) Let $f$ be a function defined on $[a,b]$. The left and right Riemann-Liouville fractional derivatives of order $\gamma>0$ for function $f$ denoted by $_aD_t^\gamma f(t)$ and $_tD_b^\gamma f(t)$, respectively, are defined by
  \begin{eqnarray*}
  &  & _aD_t^\gamma f(t)=\frac{d^{n}}{dt^{n}} {_aD_t^{\gamma-n}} f(t)=\frac{1}{\Gamma(n-\gamma)}\frac{d^{n}}{dt^{n}}\left(\int_a^t (t-s)^{n-\gamma-1} f(s)ds\right),\\
  &  & _tD_b^\gamma f(t)=(-1)^{n} \frac{d^{n}}{dt^{n}} {_tD_b^{\gamma-n}} f(t)=\frac{(-1)^{n}}{\Gamma(n-\gamma)}\frac{d^{n}}{dt^{n}}\left(\int_t^b (s-t)^{n-\gamma-1} f(s)ds\right),
 \end{eqnarray*}
 where $t\in[a,b], n-1\leq\gamma<n$ and $n\in \mathbb {N}$.
 \vskip2mm
 \noindent
 {\bf Definition 2.3.} (\cite{Jiao F2012}) Let $0<\alpha\leq1$ and $1<p<\infty$. The fractional derivative space $E_0^{\alpha,p}$ is defined by the closure of $C_0^\infty([0,T],\mathbb R^N)$ with  the norm
  $$
  \|u\|=\left(\int_0^T |_0D_t^\alpha u(t)|^pdt+\int_0^T |u(t)|^pdt\right)^{1/p}, \ \ \forall\  u\in E_0^{\alpha,p}.
  $$
  From the definition of $ E_0^{\alpha,p}$, it is apparent that the fractional derivative space $ E_0^{\alpha,p}$ is the space of functions $u:[0,T]\to\mathbb R^N$ which is absolutely continuous and has an $\alpha$-order left and right Riemann-Liouville fractional derivative $_0D_t^\alpha u\in L^p([0,T],\mathbb R^N)$ and $u(0)=u(T)=0$ and one can define the norm on $L^p([0,T],\mathbb R^N)$ as
  $$
  \|u\|_{L^p}=\left(\int_0^T|u(t)|^pdt\right)^{1/p}.
  $$
  $E_0^{\alpha,p}$ is uniformly convex by the uniform convexity of $L^p$.
 \vskip2mm
 \noindent
 {\bf Remark 2.1.} It is easy to see that $\|u\|_{V}$ defined by (\ref{P1}) is also a norm on $E_0^{\alpha,p}$
 and  $ \|u\|_{V}$ and $\|u\|$ are equivalent and
 \begin{eqnarray}\label{aaa1}
 \min\{1,V_\infty\}\|u\|^p\le \|u\|_V^p\le \max\{1,V^\infty\}\|u\|^p.
 \end{eqnarray}
 \vskip2mm
 \noindent
 {\bf Lemma 2.1.} (\cite{Jiao F2012}) Let $0<\alpha\leq1$ and $1<p<\infty$.  $E_0^{\alpha,p}$ is a reflexive and separable Banach space.
 \vskip2mm
 \noindent
 {\bf Lemma 2.2.} (\cite{Jiao F2012}) Let $0<\alpha\leq1$ and $1<p<\infty$. For all $u\in E_0^{\alpha,p}$, there has
  \begin{equation*}
   \|u\|_{L^{p}}\leq C_{p}\|_0D_t^{\alpha}u\|_{L^{p}},
   \end{equation*}
   where
   $$
   C_{p}=\frac{T^\alpha}{\Gamma(\alpha+1)}>0.
   $$
    Moreover, if $\alpha>\frac{1}{p}$, then
   \begin{equation}\label{bb2}
   \|u\|_{\infty}\leq \frac{T^{\alpha-\frac{1}{p}}}{\Gamma(\alpha)(\alpha q-q+1)^{\frac{1}{q}}}\|_0D_t^{\alpha}u\|_{L^{p}},\,\ \frac{1}{p}+\frac{1}{q}=1.
   \end{equation}
   \vskip2mm
 \noindent
 {\bf Lemma 2.3.} (\cite{Jiao F2012}) Let $1/p<\alpha\leq1$ and $1<p<\infty$. The imbedding of $E_0^{\alpha,p}$ in $C([0,T],\R^{N})$ is compact.

 \vskip2mm
 \par
Let $X$ be a Banach space. $\varphi \in C^{1}(X,\R)$ and $c\in\R$. A sequence $\{u_{n}\}\subset X$ is called (PS)$_{c}$  sequence (named after R. Palais and S. Smale) if the sequence $\{u_{n}\}$ satisfies
 \begin{equation*}
   \varphi(u_{n})\to c,\,\ \varphi'(u_{n})\to 0 .
   \end{equation*}

   \vskip2mm
 \noindent
 {\bf Lemma 2.4.} (Mountain Pass Theorem \cite{Ra1986, Willem M1997}) Let X be a Banach space, $\varphi \in C^{1}(X,\R)$, $w\in X$ and $r>0$ be such that  $\|w\|>r$ and
 $$
 b:=\inf_{\|u\|=r} \varphi(u)>\varphi(0)\geq \varphi(w).
 $$
 Then there exists a (PS)$_{c}$ sequence with
 $$
 c:=\inf_{\gamma\in\Gamma}\max_{t\in[0,1]}\varphi(\gamma(t)),
 $$
 $$
 \Gamma:=\{\gamma\in([0,1],X]):\gamma(0)=0,\gamma(1)=w\}.
 $$

\vskip2mm
 \par
 As in \cite{Chen2018}, for each $\lambda>0$,  we can define the  functional $I_{\lambda}: E_0^{\alpha,p}\to \R$ as
 $$
 I_{\lambda}(u)=\frac{1}{bp^{2}}\left(a+b\int_{0}^{T}(|_{0}D_{t}^{\alpha}u(t)|^{p}
     +V(t)|u(t)|^{p})dt\right)^{p}-\lambda \int_{0}^{T}F(t,u(t))dt-\frac{a^{p}}{bp^2}.
 $$
 It is easy to see that the assumption (H0)-(H2) can not ensure that $I_{\lambda}$ is well defined on $E_0^{\alpha,p}$. So we follow the idea  in $\cite{costa}$ and  simply sketch the outline of proof here. We use Lemma 2.4 to complete the proof. Since $F$ satisfies the growth condition only near $0$ by (H0)-(H2), in order to use the conditions globally, we modify and extend $F$ to $\bar{F}$ defined in section 3, and the corresponding functional is defined as $\bar{I}_{\lambda}$. Next we prove that $\bar{I}_{\lambda}$ has mountain pass geometry on $E_0^{\alpha,p}$. Then Lemma 2.4 implies that $\bar{I}_{\lambda}$ has a (PS)$_{c_\lambda}$ sequence. Then by a standard analysis, a convergent subsequence of the (PS)$_{c_\lambda}$ sequence is obtained to ensure that $c_\lambda$ is the critical value of  $\bar{I}_{\lambda}$. Finally, by an estimate about $\|u_\lambda\|_{\infty}$, we obtain that the critical point $u_\lambda$ of $\bar{I}_{\lambda}$ with $\|u_\lambda\|_{\infty}\leq \delta/2$ is just right  the solution of  system (\ref{aa2}) for all $\lambda>\lambda^*$ for some concrete $\lambda^*$.

 \vskip2mm
 {\section{Proofs}}
  \setcounter{equation}{0}
  \par
  Define $m(s)\in C^1(\R,[0,1])$ as an even cut-off function satisfying $sm'(s)\leq0$ and
 \begin{eqnarray}\label{cc1}
 m(s)= \begin{cases}
           1, \;\;\;  \text { if }\;\;|s| \leqslant \delta/2, \\
           0, \;\;\;  \text { if }\;\;|s| \geqslant  \delta.
          \end{cases}
   \end{eqnarray}
  Define $\bar{F}:[0,T]\times\R^N\to\R$ as
   $$
   \bar{F}(t,x)=m(|x|)F(t,x)+(1-m(|x|))M_2|x|^{q_{2}}.
   $$
   \par
   We define the  variational functional corresponding to $\bar{F}$ as
   \begin{eqnarray}\label{b1}
          \bar{I}_{\lambda}(u)
  &  =  & \frac{1}{bp^{2}}\left(a+b\int_{0}^{T}(|_{0}D_{t}^{\alpha}u(t)|^{p}
          +V(t)|u(t)|^{p})dt\right)^{p}-\lambda \int_{0}^{T}\bar{F}(t,u(t))dt-\frac{a^{p}}{bp^2}\nonumber\\
  &  =  & \frac{1}{bp^{2}}\left(a+b\|u\|_V^p\right)^{p}-\lambda \int_{0}^{T}\bar{F}(t,u(t))dt-\frac{a^{p}}{bp^2}
   \end{eqnarray}
 for all $ u\in  E_0^{\alpha,p}$. By (H0) and the definition of $\bar{F}$, it is easy to obtain that $\bar{F}$ satisfies
 \vskip1mm
 \noindent
{\it (H0)$'$  $\bar{F}(t,x)$ is continuously differentiable in $\R^N$ for a.e. $t \in [0,T]$, measurable
in $t$ for every $x\in \R^N$, and there exists
$b\in L^1([0,T];\R^+)$ such that
\begin{eqnarray*}
&  &  |\bar{F}(t,x)| \le a_0b(t)+M_2|x|^{q_2},\\
&  &  |\nabla \bar{F}(t,x)| \le (1+m_0)a_0b(t)+M_2q_2|x|^{q_2-1}+m_0M_2|x|^{q_2}
\end{eqnarray*}
for all $x\in \R^N$ and a.e. $t \in [0,T]$, $a_0=\max_{s\in [0,\delta]}a(s)$ and $m_0=\max_{s\in [\frac{\delta}{2},\delta]}|m'(s)|$.}
\vskip1mm
\noindent
 Hence, a standard argument  shows that $\bar{I}_{\lambda}\in C^1(E_0^{\alpha,p},\R)$  and
  \begin{eqnarray*}
          \langle \bar{I}_{\lambda}'(u),v\rangle
  &  =  & \left(a+b\|u\|_V^p\right)^{p-1}\left(\int_{0}^{T}[|_{0}D_{t}^{\alpha}u(t)|^{p-2}(_{0}D_{t}^{\alpha}u(t),\  _0D_{t}^{\alpha}v(t))dt+V(t)|u(t)|^{p-2}(u(t),v(t))]dt\right)\nonumber\\
  &    &-\lambda \int_{0}^{T}(\nabla\bar{F}(t,u(t)),v(t))dt
   \end{eqnarray*}
   for all $ u, v\in  E_0^{\alpha,p}$. Hence
   \begin{eqnarray*}
          \langle \bar{I}_{\lambda}'(u),u\rangle
    =   \left(a+b\|u\|_V^p\right)^{p-1}\|u\|_V^p-\lambda \int_{0}^{T}(\nabla\bar{F}(t,u(t)),u(t))dt
   \end{eqnarray*}
    for all $ u\in  E_0^{\alpha,p}$.
 \vskip2mm
 \noindent
   {\bf Lemma 3.1.} {\it Assume that (H1)-(H2) hold. Then \\
  (H1)$'$
  \begin{eqnarray*}
 0\le \bar{F}(t,x)\leq M_2 |x|^{q_{2}},  \ \ \mbox{for all }x\in \mathbb {R}^{N};
  \end{eqnarray*}
 (H2)$'$
   \begin{eqnarray*}
 0<\theta \bar{F}(t,x)\leq (\nabla\bar{F}(t,x),x),  \ \ \mbox{for all }x\in \R^N/\{0\},
  \end{eqnarray*}
  where $\theta=\min\{q_2,\beta\}$.
 }
 \vskip2mm
 \noindent
 {\bf Proof.}
 \par
 $\bullet$\ \   If $|x|\le \frac{\delta}{2}$, then by (H1), the conclusion (H1)$'$ holds;
 \par
              \quad  If $\frac{\delta}{2}<|x|\le \delta$, by (H1), we have
 $$
 0\le\bar{F}(t,x)=m(|x|)F(t,x)+(1-m(|x|))M_2|x|^{q_{2}}\le m(|x|)M_2|x|^{q_{2}}+(1-m(|x|))M_2|x|^{q_{2}}=M_2|x|^{q_{2}};
 $$
 \par
 \quad If $|x|\ge \delta$, then by the definition of $m$, we have
 $
 \bar{F}(t,x)=M_2|x|^{q_{2}}.
 $
 \par
 $\bullet$\ \  For all $x\in \R^N/\{0\}$, we have
   $$
   \nabla \bar{F}(t,x)=m'(|x|)\frac{x}{|x|}F(t,x)+m(|x|)\nabla F(t,x)+(1-m(|x|))q_{2}M_2|x|^{q_{2}-2}x-m'(|x|)\frac{x}{|x|}M_2|x|^{q_{2}}.
   $$
   Then
   \begin{eqnarray*}
          (\nabla \bar{F}(t,x),x)
   & = & |x|m'(|x|)(F(t,x)-M_2|x|^{q_{2}})+m(|x|)(\nabla F(t,x),x)+(1-m(|x|))q_{2} M_2|x|^{q_{2}}.
   \end{eqnarray*}
   and
    \begin{eqnarray*}
      \theta \bar{F}(t,x)-(\nabla \bar{F}(t,x),x)
  & = & m(|x|)(\theta F(t,x)-(\nabla F(t,x),x))+(\theta-q_{2})(1-m(|x|))M_2|x|^{q_{2}}\\
   & &   -|x|m'(|x|)(F(t,x)-M_2|x|^{q_{2}}).
    \end{eqnarray*}
   Apparently,  the conclusion holds for $0\leq|x|\leq\delta/2$ and $|x|\geq \delta$. If $\delta/2<|x|<\delta$, by using $\theta \le q_{2}$, the conclusion (H1), (H2) and the fact $sm'(s)\leq 0$ for all $s\in\R$, we can get the conclusion  (H2)$'$.
 \qed

 \vskip2mm
 \noindent
   {\bf Lemma 3.2.} {\it $\bar{I}_{\lambda}$ satisfies the mountain pass geometry for all $\lambda>\Lambda_1$, where $\Lambda_1$ is defined in (\ref{A1}).}
   \vskip2mm
   \noindent
   {\bf Proof.}
Note that $q_2>p^2>p$. By Lemma 3.1 and  (\ref{bb2}) , we have
 \begin{eqnarray}\label{422}
         \bar{I}_{\lambda}(u)
 & = &   \dfrac{1}{bp^{2}}(a+b\|u\|_V^{p})^{p}
        - \lambda \int_{0}^{T} \bar{F}(t,u(t))dt
           -\frac{a^{p}}{bp^{2}}
             \nonumber \\
 &\geq &  \frac{a^{p}}{bp^{2}}+\frac{a^{p-1}}{p^2}\|u\|_V^{p}
           - \lambda M_2\int_{0}^{T} |u(t)|^{q_{2}}dt
           -\frac{a^{p}}{bp^{2}}
             \nonumber \\
&\geq &  \frac{a^{p-1}}{p^2}\|u\|_V^{p}
           - \lambda M_2\|u\|_\infty^{q_{2}-p}\int_{0}^{T} |u(t)|^{p}dt
             \nonumber \\
 & \geq &  \frac{a^{p-1}}{p^2}\|u\|_V^{p}
           - \lambda M_2 \left(\frac{T^{\alpha-\frac{1}{p}}}{\Gamma(\alpha)(\alpha q-q+1)^{\frac{1}{q}}}\right)^{q_2-p}\|u\|_V^{q_{2}-p}\|u\|_{L^p}^p
             \nonumber \\
 &\geq &  \frac{a^{p-1}}{p^2}\|u\|_V^{p}
           - \lambda \frac{M_2}{V_\infty} \left(\frac{T^{\alpha-\frac{1}{p}}}{\Gamma(\alpha)(\alpha q-q+1)^{\frac{1}{q}}}\right)^{q_2-p}\|u\|_V^{q_{2}}.
             \nonumber
 \end{eqnarray}
We choose $\nu_{\lambda}=\left(\frac{a^{p-1}V_\infty}{2p^2\lambda  M_2\left(\frac{T^{\alpha-\frac{1}{p}}}{\Gamma(\alpha)(\alpha q-q+1)^{\frac{1}{q}}}\right)^{q_2-p}}\right)^{\frac{1}{q_2-p}}$ for any given $\lambda>0$. Then we have
 \begin{eqnarray}\label{k4}
         \bar{I}_{\lambda}(u)
 >    d_\lambda:=   \frac{a^{p-1}}{p^2}\nu_\lambda^p
           - \lambda \frac{M_2}{V_\infty} \left(\frac{T^{\alpha-\frac{1}{p}}}{\Gamma(\alpha)(\alpha q-q+1)^{\frac{1}{q}}}\right)^{q_2-p}\nu_\lambda^{q_{2}}>0,    \;\;\;          \mbox{ \text{for all }}  \;   \|u\|_V=\nu_{\lambda}.
 \end{eqnarray}
Choose
\begin{eqnarray}\label{kkk2}
e=\left(\frac{T}{\pi}\sin\frac{\pi t}{T},0,\cdots,0\right)\in E_0^{\alpha,p}.
\end{eqnarray}
Then
\begin{eqnarray}\label{aa4}
 \|e\|_{L^p}=D:=\begin{cases}
               \left(\frac{T^{p+1}}{\pi^{p+1}}\frac{2(p-1)!!}{p!!} \right)^{\frac{1}{p}},& \mbox{if } p \mbox{ is odd},\\
               \left(\frac{T^{p+1}}{\pi^{p}}\frac{(p-1)!!}{p!!}\right)^{\frac{1}{p}}, & \mbox{if } p \mbox{ is even}
               \end{cases}
\end{eqnarray}
and
\begin{eqnarray}\label{aa5}
\|_0D_t^{\alpha}e\|_{L^{p}}\le G :=\frac{T^{p+1-p\alpha}}{\Gamma^p(2-\alpha)(p+1-p\alpha)}.
\end{eqnarray}
By (\ref{bb2}),
\begin{eqnarray}\label{aa3}
 \|e\|_\infty\le \frac{T^{\alpha-\frac{1}{p}}}{\Gamma(\alpha)(\alpha q-q+1)^{\frac{1}{q}}}\|_0D_t^{\alpha}e\|_{L^{p}}\le G_0:=\frac{T^{\alpha-\frac{1}{p}}}{\Gamma(\alpha)(\alpha q-q+1)^{\frac{1}{q}}} G .
\end{eqnarray}
  Note that
$$
\Lambda_1=\max\left\{\frac{V_\infty a^{p-1}(\Gamma(\alpha)(\alpha q-q+1)^{\frac{1}{q}}G_0)^{q_2-p}}{2p^2M_2T^{(\alpha-\frac{1}{p})(q_2-p)}(\delta \min\{1,V_\infty\} D)^{q_2-p}}, \frac{\dfrac{1}{bp^{2}}(a+\frac{b\delta^p}{G_0^p}\max\{1,V^\infty\}(D^{p}+G^{p}))^{p}}{M_1 \frac{\delta^{q_1}}{G_0^{q_1}}T^{1-\frac{q_1}{p}}D^{q_1} }\right\}.
$$
Then
 \begin{eqnarray*}\label{l2}
\|\frac{\delta}{G_0}e\|_V\ge \frac{\delta \min\{1,V_\infty\}}{G_0} \|e\|_{L^p}\ge\nu_{\lambda}
 \end{eqnarray*}
 for all $\lambda>\Lambda_1$.
 By (\ref{aa3}), we have $\|\frac{\delta}{G_0}e\|_\infty\le \delta$. By the definition of $\bar{F}$ and (H1), we have $\bar{F}(t,x)=F(t,x)\ge M_1|x|^{q_1}$ for all $|x|\leq\delta/2$, and
$$
 \bar{F}(t,x)=m(|x|)F(t,x)+(1-m(x))M_2|x|^{q_2}\ge m(|x|)M_1|x|^{q_1}+(1-m(x))M_1|x|^{q_1}=M_1|x|^{q_1}
$$
for all $\frac{\delta}{2}<|x|\leq\delta$.
 Hence, by  H$\ddot{\mbox{o}}$lder inequality, we have
 \begin{eqnarray*}
         \bar{I}_{\lambda}(\frac{\delta}{G_0}e)
 & = &   \dfrac{1}{bp^{2}}(a+b\|\frac{\delta}{G_0}e\|_V^{p})^{p}
        - \lambda \int_{0}^{T} \bar{F}(t,\frac{\delta}{G_0}e(t))dt
           -\frac{a^{p}}{bp^{2}}
             \nonumber \\
 &\leq &  \dfrac{1}{bp^{2}}(a+b\|\frac{\delta}{G_0}e\|_V^{p})^{p}
           - \lambda M_1 \int_{0}^{T} |\frac{\delta}{G_0}e(t)|^{q_{1}}dt
           -\frac{a^{p}}{bp^{2}}
             \nonumber\\
&\leq &  \dfrac{1}{bp^{2}}(a+\frac{b\delta^p}{G_0^p}\max\{1,V^\infty\}\|e\|^{p})^{p}
           - \lambda M_1 \frac{\delta^{q_1}}{G_0^{q_1}}T^{1-\frac{q_1}{p}}\|e\|_{L^p}^{q_1} \nonumber\\
&\leq &  \dfrac{1}{bp^{2}}(a+\frac{b\delta^p}{G_0^p}\max\{1,V^\infty\}(D^p+G^p))^{p}
           - \lambda M_1 \frac{\delta^{q_1}}{G_0^{q_1}}T^{1-\frac{q_1}{p}}D^{q_1} \nonumber\\
& <  & 0
 \end{eqnarray*}
for all $\lambda> \Lambda_{1}$. \qed

\vskip2mm
\par
Let $w=\frac{\delta}{G_0}e$ and $\varphi=\bar{I}_\lambda$. Then for any given $\lambda> \Lambda_{1}$, Lemma 3.2 and Lemma 2.4 imply that  $\bar{I}_{\lambda}$ has a (PS)$_{c_\lambda}$ sequence $\{u_n\}:=\{u_{n,\lambda}\}$, that is, there exists a sequence $\{u_n\}$ satisfying
\begin{eqnarray}\label{k1}
\bar{I}_{\lambda}(u_n)\to c_{\lambda}, \quad \bar{I}_{\lambda}'(u_n)\to 0, \ \ \mbox{as } n\to \infty,
\end{eqnarray}
where
 \begin{eqnarray}\label{kkk1}
 c_{\lambda}:=\inf_{\gamma\in\Gamma}\max_{t\in[0,1]}\bar{I}_{\lambda}(\gamma(t)),
\end{eqnarray}
 $$
 \Gamma:=\{\gamma\in([0,1],X]):\gamma(0)=0,\gamma(1)=w\}.
 $$
\vskip2mm
\noindent
   {\bf Lemma 3.3.} {\it  The (PS)$_{c_\lambda}$ sequence $\{u_n\}$ has a convergent subsequence.}
   \vskip0mm
   \noindent
   {\bf Proof.} By virtue of Lemma 3.1, (\ref{k1}) and $\theta=\min\{q_2,\beta\}>p^{2}$, there exists a positive constant $M>0$ such that
 \begin{eqnarray}\label{k6}
          M+\|u_{n}\|_V
 &\geq &  \bar{I}_{\lambda}(u_{n})
         -\frac{1}{\theta} \langle\bar{I}^{\prime}_{\lambda}(u_{n}),u_{n}\rangle
          \nonumber \\
 & = &   (a+b\|u_{n}\|_V^{p})^{p-1}\left[\dfrac{1}{bp^{2}}(a+b\|u_{n}\|_V^{p})
         -\frac{1}{\theta}\|u_{n}\|_V^{p}\right]
          \nonumber \\
 &  &
           - \lambda \int_{0}^{T} \left[ \bar{F}(t,u_{n})
           -\frac{1}{\theta} (\nabla \bar{F}(t,u_{n}),u_{n}) \right] dt
           -\frac{a^{p}}{bp^{2}}
             \nonumber \\
 &\geq &  (a+b\|u_{n}\|_V^{p})^{p-1}\left[\dfrac{1}{bp^{2}}(a+b\|u_{n}\|_V^{p})
           -\frac{1}{\theta}\|u_{n}\|_V^{p}\right]
            -\frac{a^{p}}{bp^{2}}
             \nonumber \\
 &\geq &  a^{p-1}
           \left[\dfrac{a}{bp^{2}}+\left(\frac{1}{p^{2}}-\frac{1}{\theta}\right)\|u_{n}\|_V^{p}\right]
           -\frac{a^{p}}{bp^{2}}
             \nonumber \\
  &  =  &  a^{p-1}\left(\frac{1}{p^{2}}-\frac{1}{\theta}\right)\|u_{n}\|_V^{p}
 \end{eqnarray}
 for $n$ large enough, which shows that $\{u_{n}\}$ is bounded in $E_{0}^{\alpha,p}$ by $p>1$.
  By Lemma 2.1, we can assume that, up to a subsequence, for some $u_\lambda\in E_{0}^{\alpha,p}$,
  \begin{eqnarray}\label{k2}
  u_{n}&\rightharpoonup& u_\lambda \,\, \mbox{in} \,\, E_{0}^{\alpha,p},\\
  u_{n}&\to& u_\lambda \,\, \mbox{in}  \,\, C([0,T],\R^{N}).\nonumber
  \end{eqnarray}
  The following argument is similar to \cite {Zhang-Tang} with some modification. Since
  \begin{eqnarray}\label{cc4}
 &     & \left\langle\overline{I}_{\lambda}^{\prime}\left(u_{n}\right), u_{n}-u_\lambda\right\rangle\nonumber\\
 &  =  &\left(a+b\|u\|_V^p\right)^{p-1}\bigg(\int_{0}^{T}(|_{0}D_{t}^{\alpha}u_{n}|^{p-2}{_{0}D}_{t}^{\alpha}u_{n},
       {_{0}D}_{t}^{\alpha} (u_{n}-u_\lambda))dt\nonumber\\
 &    & +\int_{0}^{T}V(t)(|u_{n}|^{p-2}u_{n},u_{n}-u_\lambda)dt\bigg)-\lambda \int_{0}^{T}(\nabla \bar{F}(t,u_{n}),u_{n}-u_\lambda)dt,
  \end{eqnarray}
 we have
\begin{eqnarray}\label{cc5}
& &   \left\langle\overline{I}_{\lambda}^{\prime}\left(u_{n}\right)-\overline{I}_{\lambda}^{\prime}(u_\lambda),
                    u_{n}-u_\lambda\right\rangle
                     \nonumber\\
& = &   (a+b\|u_{n}\|_V^{p})^{p-1}\bigg(\int_{0}^{T}(|_{0}D_{t}^{\alpha}u_{n}|^{p-2}{_{0}D}
         _{t}^{\alpha}u_{n},
        {_{0}D}_{t}^{\alpha}(u_{n}-u_\lambda))dt
        \nonumber\\
&   &  +\int_{0}^{T}V(t)(|u_{n}|^{p-2}u_{n},u_{n}-u_\lambda)dt\bigg)-\lambda \int_{0}^{T}(\nabla
        \bar{F}(t,u_{n}),u_{n}-u_\lambda)dt\nonumber\\
&   &   - \bigg[(a+b\|u_{\lambda}\|_V^{p})^{p-1}\bigg(\int_{0}^{T}(|_{0}D_{t}^{\alpha}u_\lambda|^{p-2}{_{0}D}_{t}^{\alpha}u_\lambda,
        {_{0}D}_{t}^{\alpha} (u_{n}-u_\lambda))dt
      \nonumber\\
&   & + \int_{0}^{T}V(t)(|u_\lambda|^{p-2}u_\lambda,u_{n}-u_\lambda)dt\bigg)-\lambda \int_{0}^{T}(\nabla
        \bar{F}(t,u_\lambda),u_{n}-u_\lambda)dt\bigg]\nonumber\\
& = &   (a+b\|u_{n}\|_V^{p})^{p-1}\bigg(\|u_{n}\|_V^{p}-\int_{0}^{T}(|_{0}D_{t}^{\alpha}
         u_{n}|^{p-2}{_{0}D}_{t}^{\alpha}u_{n},
        {_{0}D}_{t}^{\alpha}u_{\lambda})dt       \nonumber\\
&   &      -\int _{0}^{T}V(t)(|u_{n}|^{p-2}u_{n},u_\lambda)dt\bigg)
        -(a+b\|u_\lambda\|_V^{p})^{p-1}\bigg(-\|u_\lambda\|_V^{p}
        \nonumber\\
&   &  +\int_{0}^{T}(|_{0}D_{t}^{\alpha}
       u_\lambda|^{p-2}{_{0}D}_{t}^{\alpha}u_\lambda,
        {_{0}D}_{t}^{\alpha}u_{n})dt
        +\int _{0}^{T}V(t)(|u_{\lambda}|^{p-2}u_{\lambda},u_{n})dt\bigg)
        \nonumber\\
&   &  -\lambda \int_{0}^{T}(\nabla \bar{F}(t,u_{n})-\nabla
        \bar{F}(t,u_\lambda),u_{n}-u_\lambda)dt\nonumber\\
&\geq &   (a+b\|u_{n}\|_V^{p})^{p-1}\|u_{n}\|_{V}^{p}+(a+b\|u_\lambda\|_V^{p})^{p-1}\|u_\lambda\|_V^{p}
       \nonumber\\
&   &  -(a+b\|u_{n}\|_V^{p})^{p-1}
        \bigg[\|_{0}D_{t}^{\alpha}u_{n}\|_{L^{p}}^{p-1}\|_{0}D_{t}^{\alpha}u_\lambda\|_{L^{p}}
        \nonumber\\
&   &   +\bigg(\int_{0}^{T}V(t)|u_{n}|^pdt\bigg)^{(p-1)/p}\bigg(\int_{0}^{T}V(t)|u_\lambda|^pdt\bigg)^{1/p}\bigg]
        \nonumber\\
&   &  -(a+b\|u_\lambda\|_V^{p})^{p-1}\bigg[\|_{0}D_{t}^{\alpha}u_\lambda\|_{L^{p}}^{p-1}\|_{0}D_{t}^{\alpha}u_{n}\|_{L^{p}}
       \nonumber\\
&   &      +\bigg(\int_{0}^{T}V(t)|u_\lambda|^pdt\bigg)^{(p-1)/p}\bigg(\int_{0}^{T}V(t)|u_{n}|^pdt\bigg)^{1/p}\bigg]
        \nonumber\\
&   &  -\lambda \int_{0}^{T}(\nabla \bar{F}(t,u_{n})-\nabla \bar{F}(t,u_{\lambda}),u_{n}-u_\lambda)dt
        \nonumber\\
&\geq &  (a+b\|u_{n}\|_V^{p})^{p-1}\|u_{n}\|_V^{p}+(a+b\|u_\lambda\|_V^{p})^{p-1}\|u_\lambda\|_V^{p}
        \nonumber\\
&   &  -(a+b\|u_{n}\|_V^{p})^{p-1}\bigg(
         \|_{0}D_{t}^{\alpha}u_{n}\|_{L^{p}}^{p}+\int_{0}^{T}V(t)|u_{n}|^pdt\bigg)^{(p-1)/p}
       \bigg(\|_{0}D_{t}^{\alpha}u_\lambda\|_{L^{p}}^{p}+\int_{0}^{T}V(t)|u_\lambda|^pdt\bigg)^{1/p}
        \nonumber\\
&   &  -(a+b\|u_\lambda\|_V^{p})^{p-1}\bigg(\|_{0}D_{t}^{\alpha}u_\lambda\|_{L^{p}}^{p}+
       \int_{0}^{T}V(t)|u_\lambda|^pdt\bigg)^{(p-1)/p}
      \bigg(\|_{0}D_{t}^{\alpha}u_{n}\|_{L^{p}}^{p}
       +\int_{0}^{T}V(t)|u_{n}|^pdt\bigg)^{1/p}
      \nonumber\\
&   &  -\lambda \int_{0}^{T}(\nabla \bar{F}(t,u_{n})-\nabla
         \bar{F}(t,u_\lambda),u_{n}-u_\lambda)dt \nonumber\\
& =  & (a+b\|u_{n}\|_V^{p})^{p-1}\|u_{n}\|_V^{p}+(a+b\|u_\lambda\|_V^{p})^{p-1}\|u_\lambda\|_V^{p}
          \nonumber\\
&    &     -(a+b\|u_{n}\|_V^{p})^{p-1}\|u_{n}\|_V^{p-1}\|u_{\lambda}\|_{V}-(a+b\|u_\lambda\|_V^{p})^{p-1}
           \|u_{n}\|_V\|u_\lambda\|_V^{p-1}
          \nonumber\\
&     &    - \lambda \int_{0}^{T}(\nabla \bar{F}(t,u_{n})-\nonumber\nabla
             \bar{F}(t,u_\lambda),u_{n}-u_\lambda)dt \nonumber\\
&  =  & (a+b\|u_{n}\|_V^{p})^{p-1}\|u_{n}\|_V^{p-1}(\|u_{n}\|_V-\|u_\lambda\|_V)
             \nonumber\\
&      &   +(a+b\|u_\lambda\|_V^{p})^{p-1}\|u_\lambda\|_V^{p-1}(\|u_\lambda\|_V-\|u_{n}\|_V)-\lambda \int_{0}^{T}(\nabla \bar{F}(t,u_{n})
            -             \nabla \bar{F}(t,u_\lambda),u_{n}-u_\lambda)dt
              \nonumber\\
& =&    \bigg((a+b\|u_{n}\|_V^{p})^{p-1}\|u_{n}\|_V^{p-1}
            - (a+b\|u_\lambda\|_V^{p})^{p-1}
             \|u_\lambda\|_V^{p-1}\bigg)(\|u_{n}\|_V-\|u_\lambda\|_V)
               \nonumber\\
&      &    -\lambda \int_{0}^{T}(\nabla \bar{F}(t,u_{n})
            -
          \nabla \bar{F}(t,u_\lambda),u_{n}-u_\lambda)dt.
    \end{eqnarray}
   Note that
    \begin{eqnarray}\label{k5}
    \lambda \int_{0}^{T}(\nabla \bar{F}(t,u_{n})-\nabla \bar{F}(t,u_\lambda),u_{n}-u_\lambda)dt\leq  \lambda \int_{0}^{T}|\nabla \bar{F}(t,u_{n})-\nabla \bar{F}(t,u_\lambda)||u_{n}-u_\lambda|dt\to 0 ,
    \end{eqnarray}
    by $u_{n}\to u_\lambda$ in $C([0,T],\R^{N})$ and $|\nabla \bar{F}(t,u_{n})-\nonumber\nabla \bar{F}(t,u_\lambda)|$ is bounded in $[0,T]$ because of (H0)$'$ and the boundedness of $\{u_n\}$ in  $E_0^{\alpha,p}$,
    and  (\ref{k1}) and (\ref{k2}) imply that
    \begin{eqnarray}\label{k3}
    \left\langle\overline{I}_{\lambda}^{\prime}\left(u_{n}\right)-\overline{I}_{\lambda}^{\prime}(u_\lambda), u_{n}-u_\lambda\right\rangle\to 0,\,\,\mbox{as}\,{n}\to \infty.
     \end{eqnarray}
   So by (\ref{cc5}), (\ref{k5}) and (\ref{k3}), we have
    $$
    \|u_{n}\|_V\to\|u_\lambda\|_V,\,\,\mbox{as }\,{n}\to \infty.
    $$
    By the uniform convexity of  $E_{0}^{\alpha,p}$ and $u_{n}\rightharpoonup u_\lambda$, it follows from the Kadec-Klee property
    (see \cite{Hewitt Edwin1965}) and (\ref{aaa1}), $u_{n}\to u_\lambda$ in $E_{0}^{\alpha,p}$. \qed

\vskip2mm
\par
 By the continuity of $\bar{I}_{\lambda}$, we obtain that $\bar{I}_{\lambda}(u)=c_{\lambda}$, where $c_{\lambda}$ is defined by (\ref{kkk1}).
 Then (\ref{k4}) implies that $c_\lambda\ge d_\lambda>0$. Hence $u_\lambda$ is a nontrivial critical point of $\overline{I}_{\lambda}$ in  $E_{0}^{\alpha,p}$ for any given $\lambda>\Lambda_1$.
\par
Next, we  show that $u_\lambda$ precisely is the nontrivial weak solution of system (\ref{aa2}) for any given $\lambda>\lambda^*$.
In order to get this, we need to make an estimate for the critical level $c_{\lambda}$. We introduce the functional $\widetilde{J}_{\lambda}: E_{0}^{\alpha,p}\to \R$ as follows
\begin{eqnarray}\label{422}
         \widetilde{J}_{\lambda}(u)
 & = &   \dfrac{1}{bp^{2}}(a+b\|u\|_V^{p})^{p}
        - \lambda M_1\int_{0}^{T} |u(t)|^{q_{1}}dt
           -\frac{a^{p}}{bp^{2}}.
             \nonumber
 \end{eqnarray}
 \vskip2mm
 \noindent
{\bf Lemma 3.4.}  For all $\lambda\ge \max\{\Lambda_1, \Lambda_2\}$,
\begin{eqnarray*}
        c_{\lambda}\leq C_*\lambda^{-\frac{p-1}{q_{1}-p}},
 \end{eqnarray*}
 where
$C_*$ is defined by (\ref{kkk3})
which is obviously independent of $\lambda$.
\vskip2mm
\noindent
{\bf Proof.}  Define $f_i:[0,\infty)\to \R$, $i=1,2$, by
\begin{eqnarray*}
 &  &      f_{1}(s)=   \dfrac{1}{bp^{2}}(a+bs^{p}\|e_1\|_V^{p})^{p}
                      - \lambda^{\frac{1}{q_{1}}}\|e_1\|_V^{p}\frac{s^{p}}{p}
                      - \frac{a^{p}}{bp^{2}},\\
 &  &      f_{2}(s)= - \lambda M_1s^{q_{1}}\int_{0}^{T}|e_1|^{q_{1}}dt
                      + \lambda^{\frac{1}{q_{1}}}\|e_1\|_V^{p}\frac{s^{p}}{p},
 \end{eqnarray*}
where $e_1=\frac{\delta}{G_0}e$ and $e$ is defined in (\ref{kkk2}). Then $f_{1}(s)+f_{2}(s)=\widetilde{J}_{\lambda}(s e_1)$. Let
$$
f_{2}^{\prime}(s)= - \lambda M_1q_{1} \|e_1\|_{L^{q_{1}}}^{q_{1}}s^{q_{1}-1}+ \lambda^{\frac{1}{q_{1}}}\|e_1\|_V^{p}s^{p-1}=0.
$$
Thus for each given $\lambda>0$, we have
$ s= \left(\dfrac{\lambda^{\frac{1}{q_{1}}}\|e_1\|_V^{p}}{\lambda M_1 q_{1} \|e_1\|_{L^{q_{1}}}^{q_{1}}}\right)^{\frac{1}{q_{1}-p}}.$
Then
$$\max_{s\geq0}f_{2}(s)=\left(
\dfrac{1}{p(M_1q_{1})^{\frac{p}{q_{1}-p}}}-\dfrac{M_1}{(M_1q_{1})^{\frac{q_{1}}{q_{1}-p}}}
\right)
\left(\dfrac{\|e_1\|_V}{\|e_1\|_{L^{q_{1}}}}\right)^{\frac{pq_{1}}{q_{1}-p}}
\lambda^{-\frac{p-1}{q_{1}-p}}.
$$
Obviously, $f_1(0)=0$ and
\begin{eqnarray*}
 f_{1}'(s)=   (a+bs^{p}\|e_1\|_V^{p})^{p-1}\|e_1\|_V^{p}s^{p-1}
                      - \lambda^{\frac{1}{q_{1}}}\|e_1\|_V^{p}s^{p-1}.
\end{eqnarray*}
So if
\begin{eqnarray*}
   \lambda
& > & \Lambda_2:=\left[a+b[\max\{1,V^\infty\}]^p\frac{\delta^{p}}{G_0^{p}}(D^p+G^p)\right]^{q_1(p-1)}\\
& = & \left(a+b[\max\{1,V^\infty\}]^p\frac{\delta^{p}}{G_0^{p}}\|e\|^p\right)^{q_1(p-1)}\\
& \ge  & \left(a+bs^p\|e_1\|_V^p\right)^{q_1(p-1)},
\end{eqnarray*}
$f_1(s)$ is decreasing on $s\in [0,1]$ and then $f_{1}(s)<0$ for all $s\in [0,1]$.
By (\ref{aa3}), we have
\begin{eqnarray}\label{l1}
\|s e_1\|_{\infty}\le \|\frac{\delta}{G_0}e\|_\infty \le  \delta
\end{eqnarray}
for all $ s \in [0,1]$.
 Then for all $\lambda>\Lambda_2$, by (H1)$'$,  (\ref{aa4}), (\ref{aa5}) and H$\ddot{\mbox{o}}$lder inequality, we have
\begin{eqnarray*}
&    &       c_{\lambda} \leq \max_{s \in[0,1]}\bar{I}_{\lambda}(s e_1)
      \leq \max_{s \in[0,1]}\widetilde{J}_{\lambda}(s e_1)
      \leq \max_{s \in[0,1]}f_{1}(s)+\max_{s\geq0}f_{2}(s)\nonumber\\
&   & \leq \max_{s\geq0}f_{2}(s)  =  \left(\frac{1}{p(M_1q_{1})^{\frac{p}{q_{1}-p}}}-\frac{M_{1}}{(M_1q_{1})^{\frac{q_{1}}{q_{1}-p}}}\right)
       \left(\dfrac{\|e_1\|_V}{\|e_1\|_{L^{q_{1}}}}\right)^{\frac{pq_{1}}{q_{1}-p}} \lambda^{-\frac{p-1}{q_{1}-p}}\nonumber\\
&   & \leq \left(\frac{1}{p(M_1q_{1})^{\frac{p}{q_{1}-p}}}-\frac{M_1}{(M_1q_{1})^{\frac{q_{1}}{q_{1}-p}}}\right)
       \left(\dfrac{[\max\{1,V^\infty\}]^{1/p}(D^p+G^p)^{1/p}}{T^{\frac{1}{q_1}-\frac{1}{p}}\|e\|_{L^{p}}}\right)^{\frac{pq_{1}}{q_{1}-p}} \lambda^{-\frac{p-1}{q_{1}-p}}\nonumber\\
&  & = C_*\lambda^{-\frac{p-1}{q_{1}-p}}.
 \end{eqnarray*}
 \qed

\vskip2mm
\noindent
{\bf Proof of Theorem 1.1.}\ \
Note that $u_{\lambda}$ is  a critical point of $\bar{I}_{\lambda}$ with critical value
$c_{\lambda}$. Since $\langle \bar{I}'(u_\lambda),u_\lambda\rangle =0$, similar to the argument in (\ref{k6}) and by Lemma 3.4, we have
\begin{eqnarray}\label{kk6}
         \|u_{\lambda}\|_V^{p}
 &\leq&  \dfrac{p^{2}\theta}{a^{p-1}(\theta-p^{2})}  \bar{I}_{\lambda}(u_{\lambda})
           \nonumber \\
 & = &   \dfrac{p^{2}\theta}{a^{p-1}(\theta-p^{2})}  c_{\lambda}
           \nonumber \\
 &\leq & \dfrac{p^{2}\theta}{a^{p-1}(\theta-p^{2})}  C_* \lambda^{-\frac{p-1}{q_{1}-p}}.
 \end{eqnarray}
Since
$$
\lambda> \Lambda_3=\left(\frac{T^{p\alpha-1}}{\left[\Gamma(\alpha)(\alpha q-q+1)^{\frac{1}{q}}\right]^p} \cdot \dfrac{p^{2}\theta C_* }{a^{p-1}(\theta-p^{2})} \cdot \frac{2^p}{\delta^p} \right)^{\frac{q_1-p}{p-1}},
$$
by (\ref{kk6}), we have
 \begin{eqnarray}\label{kk7}
 \|u_{\lambda}\|_{\infty}\leq  \frac{T^{\alpha-\frac{1}{p}}}{\Gamma(\alpha)(\alpha q-q+1)^{\frac{1}{q}}} \|u_{\lambda}\|_V\le \delta/2.
 \end{eqnarray}
So for all $\lambda>\Lambda_3$, $|u_{\lambda}(t)|\le \|u_{\lambda}\|_{\infty}\le\delta/2$ for a.e. $t\in [0,T]$ and then $\bar{F}(t,u(t))={F}(t,u(t))$ for a.e. $t\in [0,T]$. Furthermore, $\bar{I}_\lambda(u_\lambda)=I_\lambda(u_\lambda)=c_\lambda>0$ and  $ \langle \bar{I}'(u_\lambda),v\rangle =\langle {I}'(u_\lambda),v\rangle =0$ for all $v\in  E_0^{\alpha,p}$. Thus $u_\lambda$ is precisely the nontrivial weak solution of system (\ref{aa2}) when $\lambda>\lambda^*:=\max\{\Lambda_1,\Lambda_2,\Lambda_3\}$. Note that $ p>1$ and $ q_{1}>p$. By (\ref{kk6}) and (\ref{kk7}), it is obvious that
\begin{eqnarray*}
 \lim_{\lambda\to \infty} \|u_{\lambda}\|_V=0=\lim_{\lambda\to \infty} \|u_{\lambda}\|_\infty.
 \end{eqnarray*}
 \qed
\vskip2mm
 {\section{Example}}
   \setcounter{equation}{0}
  Assume that $N=2, a=b=T=1$, $p=3$ and $\delta=1$. Then $q=\frac{3}{2}$. Let $q_{1}=12, q_{2}=10, F(t,x)=(t+1)|x|^{11}$ for a.e. $t\in [0,1]$ and all $x\in \R^N$ with $|x|\le 1$.  $V(t)=7t^2+1$ for all $t\in [0,1]$. Then $V^{\infty}=8$ and $V_{\infty}=1$. Choose $\alpha=\frac{1}{2}$.  Consider the system
  \begin{eqnarray}\label{ccc1}
   \begin{cases}
 A(u(t))[_{t}D_{1}^{1/2}\phi_{3}(_{0}D_{t}^{1/2}u(t))+(7t^2+1)\phi_{3}(u(t))]
     =11\lambda(t+1)|u|^{9}u, \;\;\; \mbox{a.e. } t\in[0,1],\\
   u(0)=u(1)=0,
   \end{cases}
 \end{eqnarray}
 where
 $$
 A(u(t))= \left[1+\int_{0}^{1}(|_{0}D_{t}^{1/2}u(t)|^{3}
     +(7t^2+1)|u(t)|^{3})dt\right]^{2}.
 $$
 By Theorem 1.1, we can obtain that system (\ref{ccc1}) has at least a nontrivial solution $u_{\lambda}$ in $E_{0}^{\frac{1}{2},3}$  for each $\lambda>183.46^{24}$ and
 $\lim_{\lambda\to \infty} \|u_{\lambda}\|_V=0=\lim_{\lambda\to \infty} \|u_{\lambda}\|_\infty$.
 \par
 In fact, we can verify that $F(t,u)$ satisfies (H0)-(H2) as follows.
 \par
 i) Note that
\begin{eqnarray*}
 |F(t,x)|= (t+1)|x|^{11},\ \  |\nabla F(t,x)|=11(t+1)|x|^{10}
 \end{eqnarray*}
 for all $|x|\leq \delta$ . Set $a(|x|)=|x|^{10}+1$ for all $x\in \R^N$ with $|x|\le 1$,  $b(t)=11(t+1)$ for a.e $t\in [0,T]$. Then assumption (H0) is satisfied.
 \par
 ii) Note that $q_{1}=12>q_{2}=10>p^{2}=9 $, and
 $$|x|^{12}\leq F(t,x)=(t+1)|x|^{11}\leq 2|x|^{10},$$
 for all $|x|\leq \delta$ and a.e. $t\in [0,1]$. Set $M_{1}=1$ and  $M_{2}=2$. Then assumption (H1) is also satisfied.
 \par
 iii) Let $\beta=10>p^{2}=9$. Then
 $$
 0\leq 10(t+1)|x|^{11}=\beta F(t,x)\leq 11(t+1)|x|^{11}=(\nabla F(t,x),x)
 $$
 holds for all $x\in \R^{2}$ and a.e. $t\in [0,1]$, and so assumption (H3) is satisfied.
Next, we compute the value of $\lambda^*$ by the formulas in Theorem 1.1.
Note that $\Gamma(\frac{1}{2})=\sqrt{\pi}, \,\Gamma(2-\frac{1}{2})=\frac{\sqrt{\pi}}{2}$.
 We obtain
 \begin{eqnarray*}
 D=\left(\frac{T^{p+1}}{\pi^{p+1}}\frac{2(p-1)!!}{p!!} \right)^{\frac{1}{p}}=\bigg(\frac{4}{3}\bigg)^{\frac{1}{3}}\pi^{-\frac{4}{3}},
 \end{eqnarray*}
  \begin{eqnarray*}
G=\bigg(\frac{T^{p+1-p\alpha}}{\Gamma^p(2-\alpha)(p+1-p\alpha)}\bigg)^{1/p}=\bigg(\frac{16}{5}\bigg)^
{\frac{1}{3}}\pi^{-\frac{1}{2}} ,
\end{eqnarray*}
 \begin{eqnarray*}
G_0=\frac{T^{\alpha-\frac{1}{p}}}{\Gamma(\alpha)(\alpha q-q+1)^{\frac{1}{q}}} G=5^{-\frac{1}{3}}\cdot 16^{\frac{2}{3}}\pi^{-1}.
\end{eqnarray*}
Then by $\theta=min\{\beta,q_{2}\}=10$, (\ref{A1}), (\ref{A2}) and (\ref{A3}), we have
\begin{eqnarray*}
\Lambda_{1}=\max\left\{ \sqrt[3]{\frac{768}{78125}}\pi^{35/6}, \frac{16^{11/2}\cdot3\pi^{4}}{625}(1+\frac{5}{24\pi}+\frac{\pi^{3/2}}{2})^{3}\right\},
\end{eqnarray*}
\begin{eqnarray*}
\Lambda_{2}=\bigg(1+\frac{40}{3}\pi^{-1}+32\pi^{3/2}\bigg)^{24},
\end{eqnarray*}
\begin{eqnarray*}
\Lambda_{3}=\bigg(\frac{720\cdot16C_{*}}{\pi^{3/2}}\bigg)^{\frac{9}{2}},
\end{eqnarray*}
and by (\ref{kkk3}),
\begin{eqnarray*}
C_*=16\bigg(\frac{1}{3\cdot\sqrt[3]{12}}-\frac{1}{\sqrt[3]{12^{4}}}\bigg)\bigg(1+\frac{12}{5}\pi^{5/2}\bigg)^{4/3}.
\end{eqnarray*}
 Compared $\Lambda_{1}$, $\Lambda_{2}$ and $\Lambda_{3}$, it is easy to see $\lambda^*=\Lambda_{2}\approx183.46^{24}$.
\qed

 \vskip2mm
 \noindent
 {\bf Acknowledgement}\\
This project is supported by Yunnan Ten Thousand Talents Plan Young \& Elite Talents Project, Candidate Talents Training Fund of Yunnan Province (No: 2016PY027) and National Natural Science Foundation of China (11301235).

 {}
 \end{document}